\documentclass{bioinfo}
\copyrightyear{2009} \pubyear{2009}
\usepackage{amssymb}
\usepackage{url}
\usepackage{mathrsfs}
\usepackage{color}

\usepackage{eucal}%    caligraphic-euler fonts: \mathcal{ }
\firstpage{1}

\theoremstyle{remark}

\makeatletter \@addtoreset{equation}{section}
\title[Target prediction and a statistical sampling algorithm for RNA-RNA interaction]
      {Target prediction and a statistical sampling algorithm for RNA-RNA interaction}

\author[F.W.D.\ Huang, J.\ Qin, C.M.\ Reidys, P.F.\ Stadler]{%
  Fenix W.D. Huang\,$^1$,
  Jing Qin\,$^1$,
  Christian M.\ Reidys\,$^{1,2}$\footnote{to
    whom correspondence should be addressed.
    Phone: *86-22-2350-6800;
    Fax:   *86-22-2350-9272;
    \texttt{duck@santafe.edu}%
  }, and\newline
  Peter F. Stadler$^{3-7}$}

\address{%
  $^{1}$Center for Combinatorics, LPMC-TJKLC, Nankai University
        Tianjin 300071, P.R.~China\\
  $^{2}$College of Life Science, Nankai University
         Tianjin 300071, P.R.~China\\
  $^{3}$Bioinformatics Group, Department of Computer Science,
    and Interdisciplinary Center for \phantom{$^{1}$}Bioinformatics,
    University of Leipzig,
    H{\"a}rtelstrasse 16-18, D-04107 Leipzig, Germany.\\
  $^{4}$Max Planck Institute for Mathematics in the Sciences,
    Inselstrasse 22, D-04103 Leipzig, Germany\\
  $^{5}$RNomics Group, Fraunhofer Institut for Cell Therapy and
    Immunology, Perlickstra{\ss}e 1,D-04103 \phantom{$^{1}$}Leipzig,
    Germany\\
  $^{6}$Inst.\ f.\ Theoretical Chemistry, University of Vienna,
     W{\"a}hringerstrasse 17, A-1090 Vienna, Austria\\
  $^{7}$The Santa Fe Institute, 1399 Hyde Park Rd., Santa Fe, New Mexico, USA
}

\history{Received on *****; revised on *****; accepted on *****}
\editor{Associate Editor: *****}
\begin{document}

\maketitle

\begin{abstract}
  It has been proven that the accessibility of the target sites has a
  critical influence for miRNA and siRNA. In this paper, we present a
  program, {\tt rip2.0}, not only the energetically most favorable
  targets site based on the hybrid-probability, but
  also a statistical sampling structure to illustrate the statistical
  characterization and representation of the
  Boltzmann ensemble of RNA-RNA interaction structures. The outputs
  are retrieved via backtracing an improved dynamic programming
  solution for the partition function based on the approach of Huang
  \textit{et al.} (Bioinformatics). The $O(N^6)$ time
  and $O(N^4)$ space algorithm is implemented in C (available from
  \url{http://www.combinatorics.cn/cbpc/rip2.html}).
\end{abstract}

%\keywords{RNA-RNA interaction, joint structure, dynamic programming, partition
%function, base pairing probability, loop, RNA secondary structure.}
\maketitle
%%%
%%%%%%%%%%%%%%%%%%%%%%%%%%%%%%%%%%%%%%%%%%%%%%%%%%%%%%%%%%%%%%%%%%%%%%%%%
%%%

\section{Introduction}\label{S:Introduction}

%%%
%%%%%%%%%%%%%%%%%%biological background for target prediction%%%%%%%%%%%%%%%%%%%%%%%%%%%%
%%%
Noncoding RNAs have been found to have roles in a great variety of
processes, including transcriptional regulation, chromosome
replication, RNA processing and modiTcation, messenger RNA stability
and translation, and even protein degradation and translocation.
Direct base-pairing with target RNA or DNA molecules is central to
the function of some ncRNAs \citep{Storz:02}.  Examples include the
regulation of translation in both prokaryotes \citep{Vogel:07} and
eukaryotes \citep{McManus,Banerjee}, the targeting of chemical
modifications \citep{Bachellerie}, as well as insertion editing
\citep{Benne}, transcriptional control \citep{Kugel}. The common
theme in many RNA classes, including miRNAs, siRNAs, snRNAs, gRNAs,
and snoRNAs is the formation of RNA-RNA interaction structures that
are more complex than simple sense-antisense interactions.

The hybridization energy is a widely used criterion to predict
RNA-RNA interactions \citep{rehmsmeier:04,Tjaden:06,Busch:08}. It
has been proven that the accessibility of the target sites has a
critical influence for miRNA and siRNA \citep{Ameres:07, Kertesz:07,
Kretschmer:03}. Although a lot regulatory ncRNAs has already been
identified, the number of experimentally verified target sites is
much smaller, which stimulate a great demand to restrain the list of
putative targets.
%%%%%%%%%%%%%%%%%%rip solutions and their disadvantage%%%%%%%%%%%%%%%%%%%%%%%%%%%%%%%%%%%%%%
In its most general form, the RNA-RNA interaction problem (RIP) is
NP-complete \citep{Alkan:06,Mneimneh:07}. The argument for this
statement is based on an extension of the work of \cite{Akutsu} for
RNA folding with pseudoknots. Polynomial-time algorithms can be
derived, however, by restricting the space of allowed configurations
in ways that are similar to pseudoknot folding algorithms
\citep{Rivas}. The second major problem concerns the energy
parameters since the standard loop types (hairpins, internal and
multiloops) are insufficient; for the additional types, such as
kissing hairpins, experimental data are virtually absent. Tertiary
interactions, furthermore, are likely to have a significant impact.

Several circumscribed approaches of target prediction have been
considered in the literature. The simplest approach concatenates the
two interacting sequences and subsequently employs a slightly
modified standard secondary structure folding algorithm. For
instance, the algorithms \texttt{RNAcofold}
\citep{Hofacker,Bernhart}, \texttt{pairfold} \citep{Andronescu}, and
\texttt{NUPACK} \citep{Ren} subscribe to this strategy. The main
problem of this approach is that it cannot predict important motifs
such as kissing-hairpin loops. The paradigm of concatenation has
also been generalized to the pseudoknot folding algorithm of
\cite{Rivas}. The resulting model, however, still does not generate
all relevant interaction structures \citep{Backofen,Reidys:frame}.
An alternative line of thought is to neglect all internal
base-pairings in either strand and to compute the minimum free
energy (mfe) secondary structure for their hybridization under this
constraint. For instance, \texttt{RNAduplex} and \texttt{RNAhybrid}
\citep{rehmsmeier:04} follows this line of thought. \texttt{RNAup}
\citep{Mueckstein:05a,Mueckstein:08a} and \texttt{intaRNA}
\citep{Busch:08} restrict interactions to a single interval that
remains unpaired in the secondary structure for each partner. Due to
the highly conserved interaction motif, snoRNA/target complexes are
treated more efficiently using a specialized tool \citep{Tafer:09x}
however. \cite{Pervouchine:04} and \cite{Alkan:06} independently
derived and implemented minimum free energy (mfe) folding algorithms
for predicting the joint secondary structure of two interacting RNA
molecules with polynomial time complexity. In their model, a ``joint
structure'' means that the intramolecular structures of each
molecule are pseudoknot-free, the intermolecular binding pairs are
noncrossing and there exist no so-called ``zig-zags''. The optimal
``joint structure'' can be computed in $O(N^6)$ time and $O(N^4)$
space by means of dynamic programming.

%%%%%%%%%%%%%%%%%%rip partition function part%%%%%%%%%%%%%%%%%%%%%%%%%%%%%%%%%%%%%%

Recently, \cite{Backofen} and \cite{rip:09} independently presented
{\texttt piRNA} and {\tt rip}1.0, tools that use dynamic programming
algorithm to compute the partition function of ``joint structures'',
both in $O(N^6)$ time. Albeit differing in design details, they are
equivalent. In addition, \cite{rip:09} identified in {\tt rip}1.0 a
basic data structure that forms the basis for computing additional
important quantities such as the base pairing probability matrix.
However, since the probabilities of hybrid is not simply a sum of
the probabilities of the exterior arcs which are not independent,
{\tt rip}1.0 can not solve the probability of a hybrid.
%%%
%%%%%%%%%%%%%%%%%%%%%%%%%hybrid-matrix%%%%%%%%%%%%%%%%%%%%%%%%%%%%%%%%%%%%%%%%%%
%%%
\begin{figure}[ht]
\begin{center}
\includegraphics[width=\columnwidth]{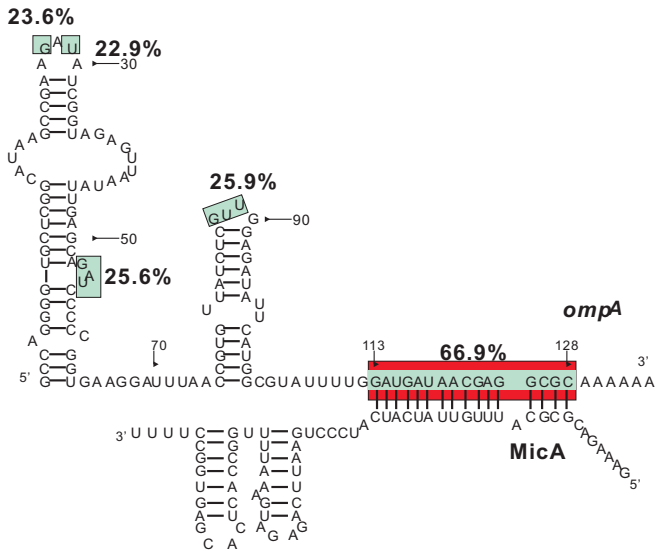}
\end{center}
\caption{\small The natural  structure of \emph{ompA}-\emph{MicA}
\citep{Udekwu:05}, in which the target site are colored in red and
the regions colored in green are the ones with the first five
region-probabilities. The target sites $R[i,j]$ of $ompA$ (interacts
with $MicA$) whose probabilities larger than $10^{-1}$ are showed in
Tab.~\ref{T:ompA}.} \label{F:ompA}
\end{figure}
%%%
%%%%%%%%%%%%%%%%%%%%%%%%%%%%%%%%%%%%%%%%%%%%%%%%%%%%%%%%%%%%%%%%%%%
%%%
%%%
%%%%%%%%%%%%%%%%%%%%%%%%%hybrid-matrix%%%%%%%%%%%%%%%%%%%%%%%%%%%%%%%%%%%%%%%%%%
%%%
\begin{figure}[ht]
\begin{center}
\includegraphics[width=0.9\columnwidth]{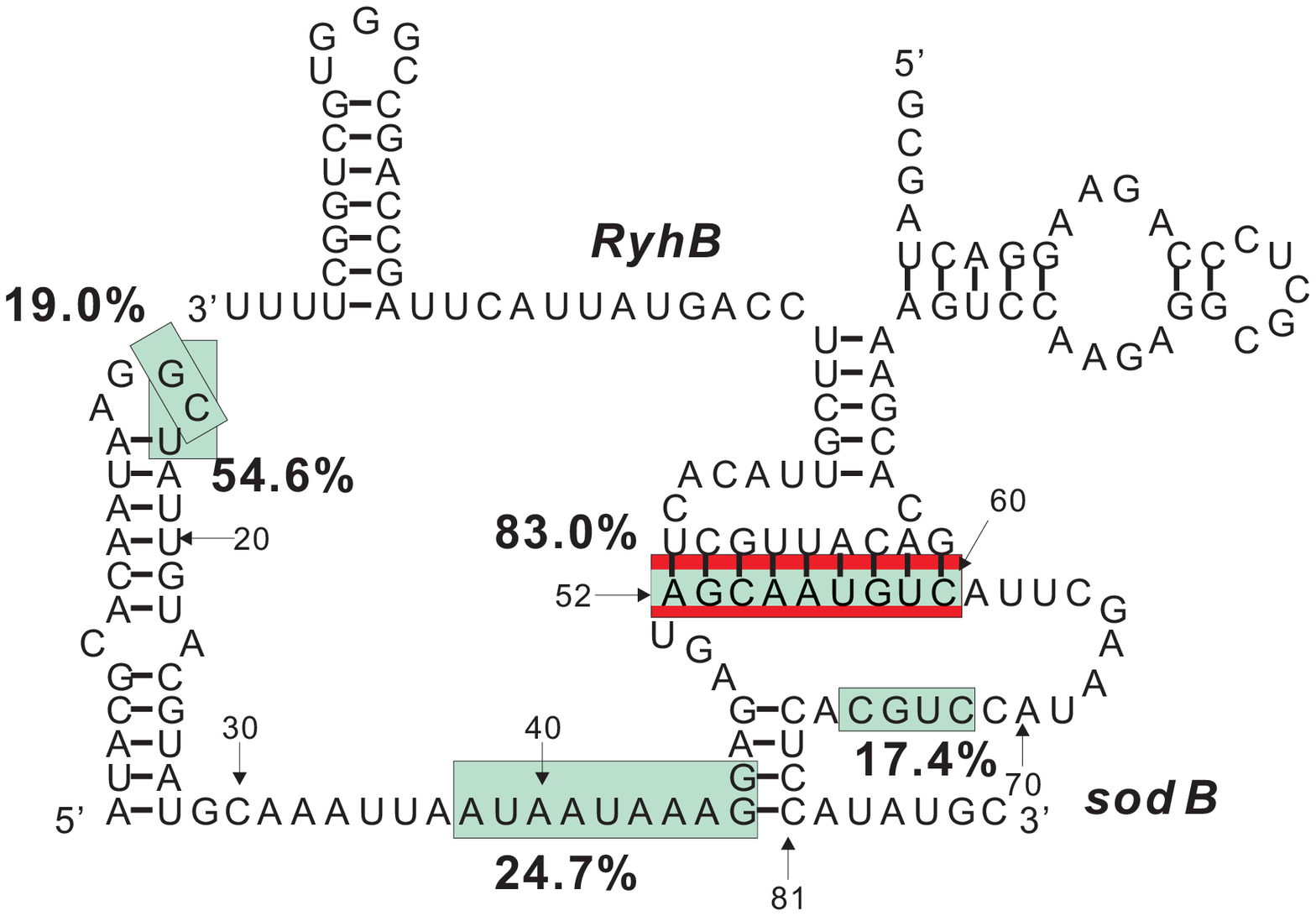}
\end{center}
\caption{\small The natural  structure of \emph{sodB}-\emph{RyhB}
\citep{Geissmann}, in which the target site are colored in red and
the regions colored in green are the ones with the first five
region-probabilities. The target sites $R[i,j]$ of \emph{sodB}
(interacts with \emph{RyhB}) whose probabilities larger than
$10^{-1}$ are showed in Tab.~\ref{T:sodB}.} \label{F:sodB}
\end{figure}
%%%
%%%%%%%%%%%%%%%%%%%%%%%%%%%%%%%%%%%%%%%%%%%%%%%%%%%%%%%%%%%%%%%%%%%
%%%
%%%
%%%%%%%%%%%%%%%%%%%%%%%%%hybrid-matrix%%%%%%%%%%%%%%%%%%%%%%%%%%%%%%%%%%%%%%%%%%
%%%
\begin{figure}[ht]
\begin{center}
\includegraphics[width=0.9\columnwidth]{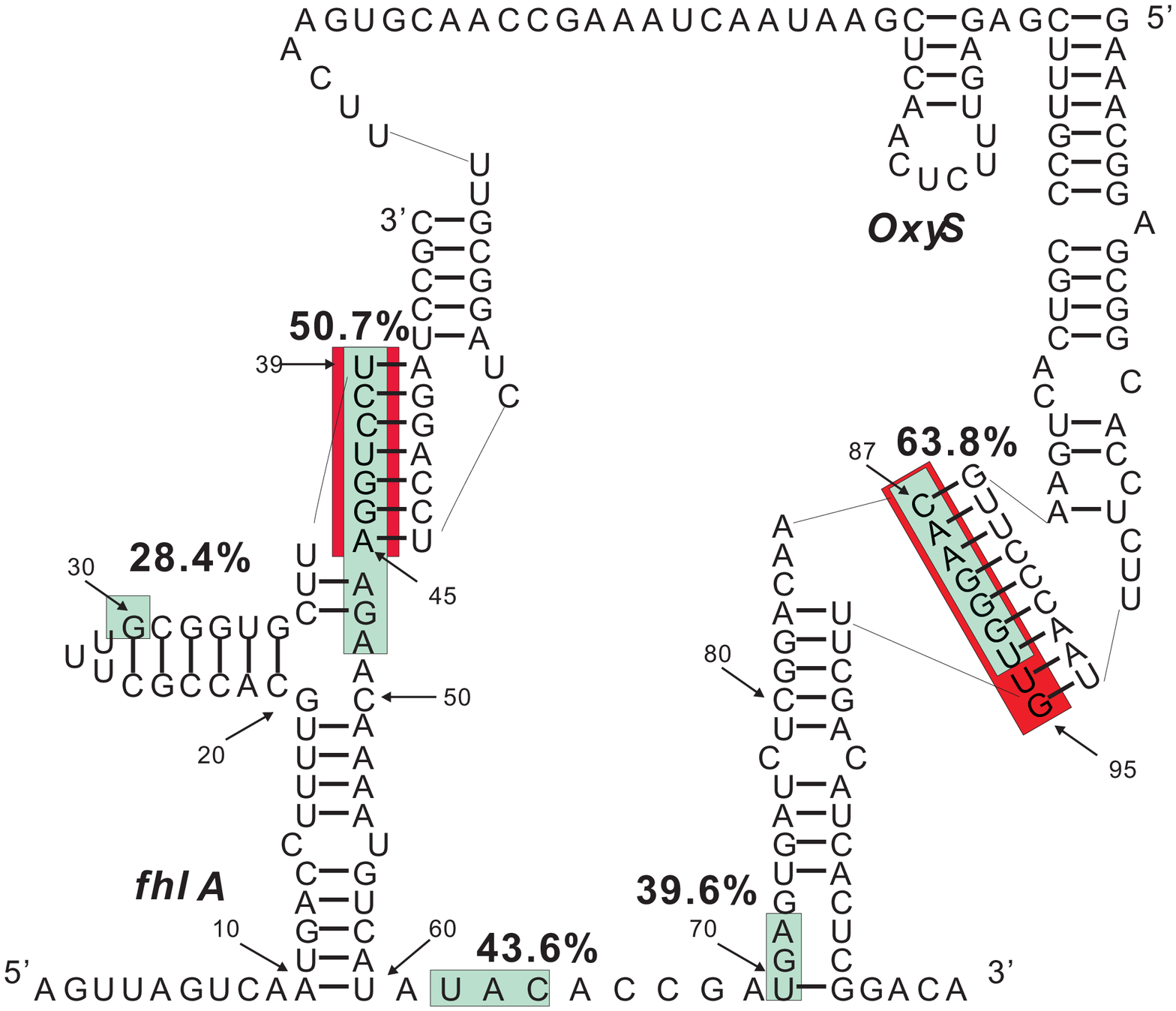}
\end{center}
\caption{\small The natural  structure of \emph{fhlA}-\emph{OxyS}
\citep{Backofen}, in which the target site are colored in red and
the regions colored in green are the ones with the first five
region-probabilities. The target sites $R[i,j]$ of \emph{fhlA}
(interacts with \emph{OxyS}) whose probabilities larger than
$10^{-1}$ are showed in Tab.~\ref{T:fhlA}.} \label{F:fhlA}
\end{figure}
%%%
%%%%%%%%%%%%%%%%%%%%%%%%%%%%%%%%%%%%%%%%%%%%%%%%%%%%%%%%%%%%%%%%%%%
%%%
%%%%%%%%%%%%%%%%%%from pf to statistic sampling%%%%%%%%%%%%%%%%%%%%%%%%%%%%%%%%%%%%%%

The calculation of equilibrium partition functions and base-pairing
probabilities is an important advance toward the characterization of
the Boltzmann ensemble of RNA-RNA interaction structures. However,
this elegant algorithm does not generate any structures. However, as
\cite{Ding:03} suggested with prototype algorithms, the generation
of a statistically representative sample of secondary structures may
provide a resolution to this dilemma.

%%%%%%%%%%%%%%%%%%rip 1.0 and 2.0 differences%%%%%%%%%%%%%%%%%%%%%%%%%%%%%%%%%%%%%%

In contrast to {\tt rip}1.0, given two RNA sequences, the output of
{\tt rip}2.0 consists of not only the partition function, the base
pairing probability matrix,  but also the contact-region probability
matrix based on the hybrid probabilities via introducing a new
component ``hybrid" in the decomposition process and a statistically
sampled RNA-RNA structure based on the probability-matrices. At the
same time, we decrease the storage space from 4D-matrices and
2D-matrices to 4D-matrices and 2D-matrices.
\begin{methods}
%%%
%%%%%%%%%%%%%%%%%%%%%%%%%%%%%%%%%%%%%%%%%%%%%%%%%%%%%%%%%%%%%%%%%%%%%
%%%
\section{Partition function}
%%%
%%%%%%%%%%%%%%%%%%%%%%%%%%%%%%%%%%%%%%%%%%%%%%%%%%%%%%%%%%%%%%%%%%%%%
%%%
\subsection{Background}
%%%
%%%%%%%%%%%%%%%%%%%%%%%%%%%%%%%%%%%%%%%%%%%%%%%%%%%%%%%%%%%%%%%%%%%%%
%%%
Let us first review some basic concepts has been introduced by
\cite{rip:09}, see supplement material (SM) for a full-version.

Given two RNA sequences $R$ and $S$ (e.g.~an antisense RNA and its
target) with $N$ and $M$ vertices, we index the vertices such that
$R_1$ is the $5'$ end of $R$ and $S_1$ denotes the $3'$ end of $S$.
The edges of $R$ and $S$ represent the intramolecular base pairs. A
\emph{joint structure}, $J(R,S,I)$, is a graph with the following
properties, see Fig.~\ref{F:joint}, \textbf{(B)}:
\begin{enumerate}
\item $R$, $S$ are secondary structures (each nucleotide being paired
with at most one other nucleotide via hydrogen bonds, without
internal pseudoknots);
\item $I$ is a set of arcs of the form $R_iS_j$ without
   pseudoknots, i.e., if $R_{i_1}S_{j_1}$, $R_{i_2}S_{j_2}\in I$
   where $i_1<i_2$, then $j_1<j_2$ holds;
\item There are no ¡°zig-zags¡±, see Fig.~\ref{F:joint},
\textbf{(A)}.
\end{enumerate}
Joint structures are exactly the configurations that are considered
in the maximum matching approach of \cite{Pervouchine:04}, in the
energy minimization algorithm of \cite{Alkan:06}, and in the
partition function approach of \cite{Backofen}. The subgraph of a
joint structure $J(R,S,I)$ induced by a pair of subsequences $\{R_i,
R_{i+1},\dots,R_j\}$ and $\{S_h, S_{h+1},\dots,S_\ell\}$ is denoted
by $J_{i,j;h,\ell}$. In particular, $J(R,S,I)=J_{1,N;1,M}$.  We say
$R_{a}R_{b} (S_{a}S_{b},R_{a}S_{b})\in J_{i,j;h,\ell}$ if and only
if $R_{a}R_{b} (S_{a}S_{b},R_{a}S_{b})$ is an edge of the graph
$J_{i,j;h,\ell}$.  Furthermore, $J_{i,j;h,\ell}\subset J_{a,b;c,d}$
if and only if $J_{i,j;h,\ell}$ is a subgraph of $J_{a,b;c,d}$
induced by $\{R_i,\dots,R_j\}$ and $\{S_h,\dots,S_\ell\}$.
%%%
%%%%%%%%%%%%%%%%%%%%%%%%%%zigjoint%%%%%%%%%%%%%%%%%%%%%%%%%%%%%%%%%%%%%%%%%%
%%%
\begin{figure}[ht]
 \begin{center}
  \includegraphics[width=\columnwidth]{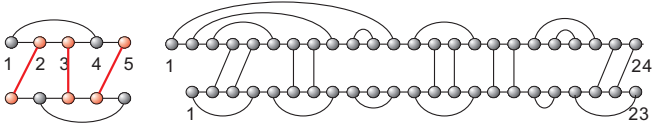}
\end{center}
\caption{ \textsf{(A):} A zigzag, generated by
  $R_2S_1$, $R_3S_3$ and $R_5S_4$ (red). \textsf{(B):}
  the joint structure $J_{1,24;1,23}$, we color the different
  segments and tight structures in which $J_{1,24;1,23}$ decomposes.
} \label{F:joint}
\end{figure}
%%%
%%%%%%%%%%%%%%%%%%%%%%%%%%%%%%%%%%%%%%%%%%%%%%%%%%%%%%%%%%%%%%%%%%%%
%%%
Given a joint structure, $J_{a,b;c,d}$, its tight structure (ts)
$J_{a',b';c',d'}$ is either a single exterior arc $R_{a'}S_{c'}$ (in
the case  $a'=b'$ and $c'=d'$), or the minimal block centered around
the leftmost and rightmost exterior arcs $\alpha_l,\alpha_r$,
(possibly being equal) and an interior arc subsuming both, i.e.,
$J_{a',b';c',d'}$ is tight in $J_{a,b;c,d}$ if it has either an arc
$R_{a'}R_{b'}$ or $S_{c'}S_{d'}$ if $a'\neq b'$ or $c'\neq d'$.

In the following, a ts is denoted by $J^{T}_{i,j;h,\ell}$. If
$J_{a',b';c',d'}$ is tight in $J_{a,b;c,d}$, then we call
$J_{a,b;c,d}$ its envelope. By construction, the notion of ts is
depending on its envelope. There are only four basic types of ts,
see Fig.~\ref{F:typeins}:
\begin{itemize}
\item[$\circ:$]        $\{R_iS_h\}=J^{\circ}_{i,j;h,\ell}$ and $i=j$,
        $h=\ell$;
\item[$\bigtriangledown:$] $R_iR_j\in J^{\bigtriangledown}_{i,j;h,\ell}$ and
      $S_{h}S_{\ell}\not\in J^{\bigtriangledown}_{i,j;h,\ell}$;
\item[$\square:$]      $\{R_iR_j,S_{h}S_{\ell}\}\in
J^{\square}_{i,j;h,\ell}$;
\item[$\bigtriangleup:$] $S_{h}S_{\ell} \in J^{\bigtriangleup}_{i,j;h,\ell}$
      and $R_iR_j\not \in J^{\bigtriangleup}_{i,j;h,\ell}$.
\end{itemize}

%%%
%%%%%%%%%%%%%%%%%%%%%%%%%%%tight-example%%%%%%%%%%%%%%%%%%%%%%%%%%%%%%%%%%%%%%%%%%%%%%
%%%
\begin{figure}[ht]
\begin{center}
  \includegraphics[width=\columnwidth]{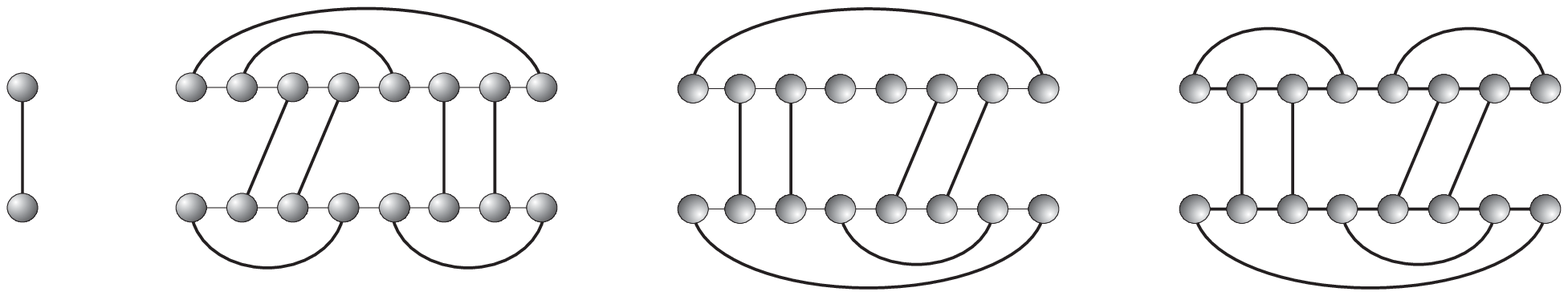}
\end{center}
\caption{From left to right: tights of type $\circ$,
$\bigtriangledown$, $\square$ and $\bigtriangleup$.}
\label{F:typeins}
\end{figure}

%%%
%%%%%%%%%%%%%%%%%%%%%%%%%%%%%%%%%%%%%%%%%%%%%%%%%%%%%%%%%%%%%%%%%
%%%
\subsection{Refined decomposition grammar for target prediction}
%%%
%%%%%%%%%%%%%%%%%%%%%%%%%%%%%%%%%%%%%%%%%%%%%%%%%%%%%%%%%%%%%%%%%
%%%

The unique ts decomposition would in principle already suffice to
construct a partition function algorithm. Indeed, each decomposition
step corresponds to a multiplicative recursion relation for the
partition functions associated with the joint structures. From a
practical point of view, however, this would result in an unwieldy
expensive implementation. The reason are the multiple break points
$a$, $b$, $c$, $d$, \dots, each of which correspond to a nested
\texttt{for}-loop.

We therefore need a refined decomposition that reduced the number of
break points. To this end we call a joint structure
\emph{right-tight} (rts), $J^{RT}_{i,j;r,s}$  in
$J_{i_1,j_1;r_1,s_1}$ if its rightmost block is a
$J_{i_1,j_1;r_1,s_1}$-ts and \emph{double-tight} (dts),
$J^{DT}_{i,j;r,s}$ in $J_{i_1,j_1;r_1,s_1}$ if both of its left- and
rightmost blocks are $J_{i_1,j_1;r_1,s_1}$-ts's. In particular, for
the convenient of the computation, we assume the single interaction
arc as a special case of dts,
i.e.~$J^{DT}_{i_1,i_1;r_1,r_1}=R_{i_1}S_{r_1}$ . In order to obtain
the probability of a hybrid $J^{\sf Hy}_{i_1,i_\ell;j_1,j_\ell}$ via
the backtracing method used in \citep{rip:09}, we introduce the
hybrid structure, $J^{\sf Hy}_{i_1,i_\ell;j_1,j_\ell}$, as a new
block item used in the decomposition process. We adopt the point of
view of Algebraic Dynamic Programming \citep{Giegerich:02} and
regard each decomposition rule as a production in a suitable
grammar. Fig.~\ref{F:grammar} summarizes three major steps in the
decomposition: (I) ``interior arc-removal'' to reduce ts. The scheme
is complemented by the usual loop decomposition of secondary
structures, and (II) ``block-decomposition'' to split a joint
structure into two blocks.

\begin{figure}[ht]
\begin{center}
\includegraphics[width=\columnwidth]{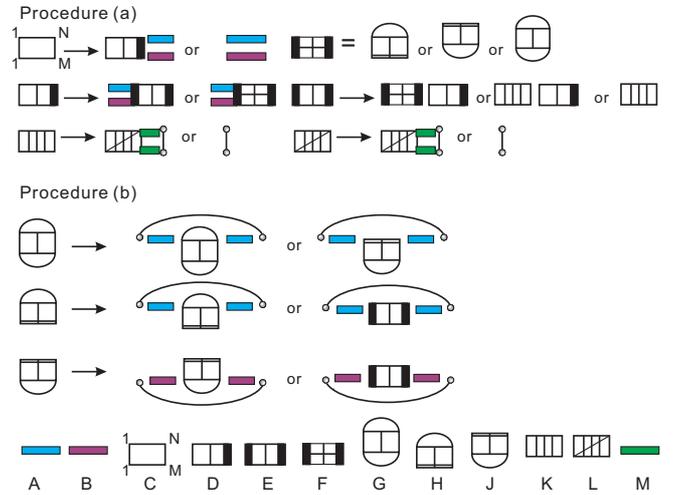}
\end{center}
\caption{Illustration of Procedure (a) the reduction of arbitrary
joint structures
  and right-tight structures, and Procedure (b)
  the decomposition of tight structures. The panel below indicates the 10
  different types of structural components: \textbf{A}, \textbf{B}:
  maximal secondary structure segments $R[i,j]$, $S[r,s]$;
  \textbf{C}: arbitrary joint structure $J_{1,N;1,M}$;
  \textbf{D}: right-tight structures $J^{RT}_{i,j;r,s}$;
  \textbf{E}: double-tight structure $J^{DT}_{i,j;r,s}$;
  \textbf{F}  tight structure of type $\bigtriangledown$,
    $\bigtriangleup$ or $\square$;
  \textbf{G}: type $\square$ tight structure $J^{\square}_{i,j;r,s}$;
  \textbf{H}: type $\bigtriangledown$ tight structure
    $J^{\bigtriangledown}_{i,j;r,s}$;
  \textbf{J}: type $\bigtriangleup$ tight structure
    $J^{\bigtriangleup}_{i,j;r,s}$;
  \textbf{K}: hybrid structure $J^{\sf Hy}_{i,j;h,\ell}$;
  \textbf{L}: substructure of a hybrid $J^{\sf
  h}_{i,j;h,\ell}$ such that $R_{i}S_{\ell}$ and $R_{h}S_{\ell}$ are
  exterior arcs;
  \textbf{M}: isolated segment $R[i,j]$ or $S[h,\ell]$.
   } \label{F:grammar}
\end{figure}

According to the decomposition rule, a given joint structure
decomposed into interior arcs and hybrids, see Figure \ref{F:tree}
\textbf{(A)}. The details of the decomposition procedures are
collected in SM, Section~2, where we show that for each joint
structure $J_{1,N;1,M}$ we indeed obtain a unique decomposition-tree
(parse-tree), denoted by $T_{J_{1,N;1,M}}$. More precisely,
$T_{J_{1,N;1,M}}$ has root $J_{1,N;1,M}$ and all other vertices
correspond to a specific substructure of $J_{1,N;1,M}$ obtained by
the successive application of the decomposition steps of
Fig.~\ref{F:grammar} and the loop decomposition of the secondary
structures. The decomposition trees of a concrete example generated
according to {\tt rip}2.0 and {\tt rip}1.0  is shown in
Fig.~\ref{F:tree} \textbf{(A)} and \textbf{(B)}, respectively.

%%%
%%%%%%%%%%%%%%%%%%%%%%%%%%%%%%%%%%%%%%%%%%%%%%%%%%%%%%%%%%%%%%%%
%%%
\begin{figure}[tb]
\begin{center}
  \includegraphics[width=\columnwidth]{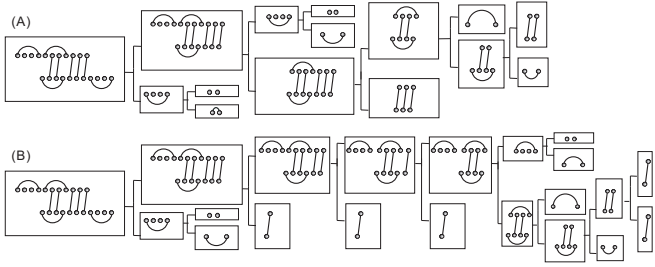}
\end{center}
\caption{The decomposition trees $T_{J_{1,15;1,8}}$ for the joint
structure $J_{1,15;1,8}$ according to the grammar in {\tt rip}2.0
\textbf{(A)} and {\tt rip}1.0 \textbf{(B)}, respectively.}
\label{F:tree}
\end{figure}
%%%
%%%%%%%%%%%%%%%%%%%%%%%%%%%%%%%%%%%%%%%%%%%%%%%%%%%%%%%%%%%%%%%%
%%%

Let us now have a closer look at the energy evaluation of
$J_{i,j;h,\ell}$. Each decomposition step in Fig.~\ref{F:grammar}
results in substructures whose energies we assume to contribute
additively and generalized loops that need to be evaluated directly.
There are the following two scenarios:

\par\noindent\textbf{I.\ Interior Arc removal.} The first type of
decomposition is focus on decomposing ts which is similar as the
approach deduced by \cite{rip:09}. Most of the decomposition
operations in Procedure (b) displayed in Fig.~\ref{F:grammar} can be
viewed as the ``removal'' of an arc (corresponding to the closing
pair of a loop in secondary structure folding) followed by a
decomposition. Both: the loop-type and the subsequent possible
decomposition steps depend on the newly exposed structural elements.
W.l.o.g., we may assume that we open an interior base pair $R_iR_j$.

For instance, a rts $J^{RT}_{p,q,r,s}$ (denoted by ``\textsf{D}'' in
Fig.~\ref{F:grammar}) we need to determine the type of the exposed
pairs of both $R[p,q]$ and $S[r,s]$. Hence each such structure will
be indexed by two types lies in $\{{\sf E},{\sf M},{\sf K},{\sf
F}\}$. Analogously, there are in total four types of a hybrid
$J^{\sf Hy}_{i,j;h,\ell}$, i.e.~ $\{J^{\sf Hy,EE}_{i,j;h,\ell},
J^{\sf Hy,EK}_{i,j;h,\ell},J^{\sf Hy,KE}_{i,j;h,\ell}, J^{\sf
Hy,KK}_{i,j;h,\ell}\}$.

%%%
%%%%%%%%%%%%%%%%%%%%%%%%%%%%%%%%%%%%%%%%%%%%%%%%%%%%%%%%%%%%%%%%%%%%
%%%
\begin{figure}[ht]
\begin{center}
\includegraphics[width=\columnwidth]{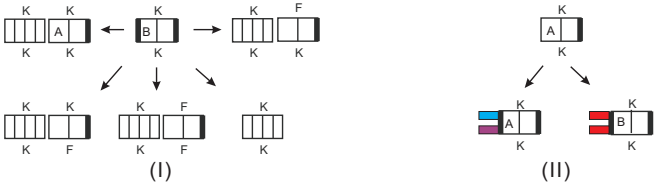}
\end{center}
\caption{({\bf I}) Decomposition of $J^{RT,{\sf MK}}_{i,j;h,\ell}$
and ({\bf II}) decomposition of $J^{DT,{\sf KKB}}_{i,j;h,\ell}$.}
\label{F:parexp3}
\end{figure}
%%%
%%%%%%%%%%%%%%%%%%%%%%%%%%%%%%%%%%%%%%%%%%%%%%%%%%%%%%%%%%%%%%%%%%%%
%%%

\par\noindent\textbf{II.\ Block decomposition.}The second type of
decomposition is the splitting of joint structures into ``blocks''.
There are two major differences in contrast to the method used in
\cite{rip:09}. First, we introduce the hybrid itself as a new block
item in the grammar and furthermore decompose a hybrid via
simultaneously removing a single exterior arc. Second, we split the
whole interaction structure into blocks via the alternating
decompositions of a rts and a dts as showed in the Procedure (a) of
Fig.~\ref{F:grammar}.

In order to make sure the maximality of a hybrid, the rts's
$J^{RT,KK}_{i,j;h,\ell}$, $J^{RT,KE}_{i,j;h,\ell}$,
$J^{RT,EK}_{i,j;h,\ell}$ and $J^{RT,EE}_{i,j;h,\ell}$ may appear in
two ways, depending on whether or not there exists an exterior arc
$R_{i_1}S_{j_1}$ such that $R[i,i_1-1]$ and $S[j,j_1-1]$ are
isolated segments. If there exists, we say rts is of type ({\sf B})
or ({\sf A}), otherwise. Similarly, a dts, $J^{RT,KK}_{i,j;h,\ell}$,
$J^{RT,KE}_{i,j;h,\ell}$, $J^{RT,EK}_{i,j;h,\ell}$ or
$J^{RT,EE}_{i,j;h,\ell}$ is of type ({\sf B}) or ({\sf A}) depending
on whether $R_iS_h$ is an exterior arc. For instance,
Fig.~\ref{F:parexp3} ({\bf I}) displays the decomposition of
$J^{DT,{\sf KKB}}_{i,j;h,\ell}$ into hybrid and rts with type ({\sf
A}) and furthermore Fig.~\ref{F:parexp3} ({\bf II}) displays the
decomposition of $J^{RT,{\sf KKA}}_{i,j;h,\ell}$.

Suppose $J^{DT}_{i,j;r,\ell}$ is a dts contained in a kissing loop,
that is we have either $E^{e}_{R[i,j]}\neq \varnothing$ or
$E^{e}_{S[h,\ell]}\neq \varnothing$. W.l.o.g., we may assume
$E^{e}_{R[i,j]}\neq \varnothing$. Then at least one of the two
``blocks" contains the exterior arc belonging to $E^{e}_{R[i,j]}$
labeled by ${\sf K}$ and ${\sf F}$, otherwise, see
Fig.~\ref{F:parexp3} ({\bf I}).

%%%
%%%%%%%%%%%%%%%%%%%%%%%%%%%%%%%%%%%%%%%%%%%%%%%%%%%%%%%%%%%%%%%%%
%%%
\subsection{Examples for partition function recursions}
%%%
%%%%%%%%%%%%%%%%%%%%%%%%%%%%%%%%%%%%%%%%%%%%%%%%%%%%%%%%%%%%%%%%%
%%%

The computation of the partition function proceeds ``from the inside
to the outside'', see equs.~(\ref{E:DDDD}). The recursions are
initialized with the energies of individual external base pairs and
empty secondary structures on subsequences of length up to four. In
order to differentiate multi- and kissing-loop contributions, we
introduce the partition functions $Q^{\sf{m}}_{i,j}$ and
$Q^{\sf{k}}_{i,j}$. Here, $Q^{\sf{m}}_{i,j}$ denotes the partition
function of secondary structures on $R[i,j]$ or $S[i,j]$ having at
least one arc contained in a multi-loop. Similarly,
$Q^{\sf{k}}_{i,j}$ denotes the partition function of secondary
structures on $R[i,j]$ or $S[i,j]$ in which at least one arc is
contained in a kissing loop. Let
$\mathbb{J}^{\xi,Y_1Y_2Y_3}_{i,j;h,\ell}$ be the set of
substructures $J_{i,j;h,\ell}\subset J_{1,N;1,M}$ such that
$J_{i,j;h,\ell}$ appears in $T_{J_{1,N;1,M}}$ as an interaction
structure of type $\xi\in\{DT, RT, \bigtriangledown,
\bigtriangleup,\square,\circ\}$ with loop-subtypes
$Y_1,Y_2\in\{\sf{M},\sf{K},\sf{F}\}$ on the sub-intervals $R[i,j]$
and $S[h,\ell]$, $Y_{3}\in\{{\sf A},{\sf B}\}$. Let
$Q^{\xi,Y_1Y_2Y_3}_{i,j;h,\ell}$ denote the partition function of
$\mathbb{J}^{\xi,Y_1Y_2Y_3}_{i,j;h,\ell}$.

For instance, the recursion for $Q^{DT, {\sf KKB}}_{i,j;h,\ell}$
displayed in Figure \ref{F:parexp3} ({\bf I}) is equivalent to:
\begin{equation}\label{E:RTMK}
\begin{split}
&Q^{RT, {\sf MK}}_{i,j;h,\ell}=\sum_{i_1,h_1}Q^{\sf
Hy,KK}_{i,i_1;h,h_1}Q^{RT, {\sf KKA}}_{i_1+1,j;h_1+1,\ell}+Q^{\sf
Hy,KK}_{i,i_1;h,h_1}Q^{RT, {\sf KF}}_{i_1+1,j;h_1+1,\ell}\\
&+Q^{\sf Hy,KK}_{i,i_1;h,h_1}Q^{RT, {\sf FF}}_{i_1+1,j;h_1+1,\ell}+
Q^{\sf Hy,KK}_{i,i_1;h,h_1}Q^{RT, {\sf
FK}}_{i_1+1,j;h_1+1,\ell}+Q^{{\sf Hy KK}}_{i,j;h,\ell}.
\end{split}
\end{equation}

In which, the recursions for $J^{\sf Hy,EE}_{i,j;h,\ell}, J^{\sf
Hy,EK}_{i,j;h,\ell}$, $J^{\sf Hy,KE}_{i,j;h,\ell}$, and $J^{\sf
Hy,KK}_{i,j;h,\ell}$ read:
\begin{equation}
  \begin{split}
Q^{\sf Hy,EE}_{i,j;h,\ell}&=\sum_{i_1,h_1} Q^{\sf
Hy,EE}_{i,i_1;h,h_1}e^{-(\sigma_{0}+\sigma G^{\sf Int}_
{i_1,h_1,j,\ell})};\\
Q^{\sf Hy,EK}_{i,j;h,\ell}&=\sum_{i_1,h_1} Q^{\sf
Hy,EK}_{i,i_1;h,h_1}e^{-(\sigma_{0}+\sigma G^{\sf Int}_
{i_1,h_1,j,\ell}+(\ell-h_1-1)\beta_3)};\\
Q^{\sf Hy,KE}_{i,j;h,\ell}&=\sum_{i_1,h_1} Q^{\sf
Hy,KE}_{i,i_1;h,h_1}e^{-(\sigma_{0}+\sigma G^{\sf Int}_
{i_1,h_1,j,\ell}+(j-i_1-1)\beta_3)};\\
Q^{\sf Hy,KK}_{i,j;h,\ell}&=\sum_{i_1,h_1} Q^{\sf
Hy,KK}_{i,i_1;h,h_1}e^{-(\sigma_{0}+\sigma G^{\sf Int}_
{i_1,h_1,j,\ell}+(j+\ell-i_1-h_1-2)\beta_3)}.
  \end{split}
  \label{E:DDDD}
\end{equation}
%%%
%%%%%%%%%%%%%%%%%%%%%%%%%%%%%%%%%%%%%%%%%%%%%%%%%%%%%%%%%%%%%%%%%%5
%%%
\section{Backtracing}\label{S:back}
%%%
%%%%%%%%%%%%%%%%%%%%%%%%%%%%%%%%%%%%%%%%%%%%%%%%%%%%%%%%%%%%%%%%%%5
%%%
\subsection{target prediction}
%%%
%%%%%%%%%%%%%%%%%%%%%%%%%%%%%%%%%%%%%%%%%%%%%%%%%%%%%%%%%%%%%%%%%%5
%%%

Given two RNA sequences, our sample space is the ensemble of all the
possible joint interaction structures. Let $Q^{I}$ denote the
partition function which sums over all the possible joint
structures. The probability measure of a given joint structure
$J_{1,N;1,M}$ is given by
\begin{equation}
\mathbb{P}_{J_{1,N;1,M}}=\frac{Q_{J_{1,N;1,M}}}{Q^{I}}.
\end{equation}

In contrast to the computation of the partition function ``from the
inside to the outside'', the computation of the
substructure-probabilities are obtained ``from the outside to the
inside'' via total probability formula (TPF). That is, the
longest-range substructures are computed first. This is analogous to
McCaskill's algorithm for secondary structures \citep{McCaskill}.

Set $J=J_{1,N;1,M}$, $T=T_{J_{1,N;1,M}}$ and let
$\Lambda_{J_{i,j;h,\ell}}=\{J\vert J_{i,j;h,\ell} \in T\}$ denote
the set of all joint structures $J$ such that $J_{i,j;h,\ell}$ is a
vertex in the decomposition tree $T$. Then
\begin{equation}
\mathbb{P}_{J_{i,j;h,\ell}}=\sum_{J\in
\Lambda_{i,j;h,\ell}}\mathbb{P}_{J}.
\end{equation}
By virtue of TPF, set $\theta_{s}$ be the possible parent-structure
of $J_{i,j;h\ell}$ and $\mathbb{P}_{J_{i,j;h,\ell}\vert \theta_{i}}$
be the conditional probability, we have
$\mathbb{P}_{J_{i,j;h,\ell}}=\sum_{s}
\mathbb{P}_{J_{i,j;h,\ell}\vert \theta_{s}}\mathbb{P}_{\theta_{s}}$.

Let $\mathbb{P}^{\xi,Y_1Y_2Y_3}_{i,j;h,\ell}$ be the probability of
$\mathbb{J}^{\xi,Y_1Y_2Y_3}_{i,j;h,\ell}$. For instance,
$\mathbb{P}^{RT,{\sf M}{\sf K}{\sf A}}_{i,j;h,\ell}$ is the sum over
all the probabilities of substructures $J_{i,j;h,\ell}\in
T_{J_{1,N;1,M}}$ such that $J_{i,j;h,\ell}\in \mathbb{J}^{RT,{\sf
M}{\sf K}{\sf A}}_{i,j;h,\ell}$, i.e.~a rts of type ${\sf A}$ and
$R[i,j]$, $S[h,\ell]$ are respectively enclosed by a multi-loop and
kissing loop. Given a component,
$\mathbb{J}^{\xi,Y_1Y_2Y_3}_{i,j;h,\ell}$ (showed in Figure
\ref{F:grammar}), we say another component
$\tilde{\mathbb{J}}^{\xi,Y_1Y_2Y_3}_{i,j;h,\ell}$ is its
\textit{parent-component} if and only if as a substructure,
$\tilde{J}^{\xi,Y_1Y_2Y_3}_{i,j;h,\ell}$ could be a parent structure
of ${J}^{\xi,Y_1Y_2Y_3}_{i,j;h,\ell}$ in the decomposition tree.
Accordingly, we say $\mathbb{J}^{\xi,Y_1Y_2Y_3}_{i,j;h,\ell}$ is the
\textit{child-component} of
$\tilde{\mathbb{J}}^{\xi,Y_1Y_2Y_3}_{i,j;h,\ell}$.  For instance,
$\mathbb{J}^{RT,KKB}_{i_1,j_1;h_1,\ell_1}$ is one of the
parent-component of $\mathbb{J}^{\sf Hy,KKB}_{i,j;h,\ell}$, see
Figure.~\ref{F:parexp3} ({\bf I}). Set $\Theta_{s}$ be one of the
possible parent-component of
$\mathbb{J}^{\xi,Y_1Y_2Y_3}_{i,j;h,\ell}$.  Accordingly, we have
\begin{equation}\label{E:Jprob2}
\mathbb{P}^{\xi,Y_1Y_2Y_3}_{i,j;h,\ell}=\sum_{s}
\mathbb{P}^{\xi,Y_1Y_2Y_3}_{i,j;h,\ell\vert
\Theta_{s}}\mathbb{P}_{\Theta_{s}},
\end{equation}
where by definition
\begin{equation}\label{E:Jprob1}
\mathbb{P}^{\xi,Y_1Y_2Y_3}_{i,j;h,\ell}=\sum_{J_{i,j;h,\ell}\in
\mathbb{J}^{\xi,Y_1Y_2Y_3}_{i,j;h,\ell}}\mathbb{P}_{J_{i,j;h,\ell}}.
\end{equation}
Furthermore, in the programme, we calculate
$\mathbb{P}^{\xi,Y_1Y_2Y_3}_{i,j;h,\ell\vert
\Theta_{s}}\mathbb{P}_{\Theta_{s}}$ for all $s$ during the
decomposition of $\Theta_{s}$. I.e.~given $\mathbb{P}_{\Theta_{s}}$,
we have
$\mathbb{P}_{\Theta_{s}}=\sum_{i}\mathbb{P}_{\epsilon_{i}\vert
\Theta_{s}}\mathbb{P}_{\Theta_{s}}$, where
$\mathbb{P}_{\epsilon_{i}\vert \Theta_{s}}$ denotes the conditional
probability of the event that given $\Theta_{s}$ is the
parent-component, $\epsilon_{i}$ is its child-component. In
particular, we have $\mathbb{P}^{\xi,Y_1Y_2Y_3}_{i,j;h,\ell\vert
\Theta_{s}}=\mathbb{P}_{\epsilon_{m}\vert \Theta_{s}}$ for some $m$.
%%%
%%%%%%%%%%%%%%%%%%%%%%%%%%%%%%%%%%%%%%%%%%%%%%%%%%%%%%%%%%%%%%%%%
%%%

Since the four subclasses of $J^{\sf Hy}_{i,j;h,\ell}$, i.e.~$J^{\sf
Hy,EE}_{i,j;h,\ell}, J^{\sf Hy,EK}_{i,j;h,\ell},J^{\sf
Hy,KE}_{i,j;h,\ell}$, and $J^{\sf Hy,KK}_{i,j;h,\ell}$ are
independent, we obtain
\begin{equation}
\mathbb{P}^{\sf Hy}_{i,j;h,\ell}=\mathbb{P}^{\sf
Hy,EE}_{i,j;h,\ell}+\mathbb{P}^{\sf
Hy,EK}_{i,j;h,\ell}+\mathbb{P}^{\sf
Hy,KE}_{i,j;h,\ell}+\mathbb{P}^{\sf Hy,KK}_{i,j;h,\ell}.
\end{equation}

Given a hybrid $J^{\sf Hy}_{i,j;h,\ell}$, recall the definition
target sites are $R[i,j]$ and $S[h,\ell]$. The probability of a
target site $R[i,j]$ is defined by
\begin{equation}
\mathbb{P}^{\sf tar}_{R[i,j]}=\sum_{h,\ell}\mathbb{P}^{\sf
Hy}_{i,j;h,\ell}.
\end{equation}
Analogously, we define $\mathbb{P}^{\sf tar}_{R[i,j]}$. We predict
the optimal interaction region with maximal probability, i.e.~
\begin{equation}
\mathbb{P}^{\sf opt}=\mathsf{max}_{i,j}\mathbb{P}^{\sf
tar}_{R[i,j]}.
\end{equation}
%%%
%%%%%%%%%%%%%%%%%%%%%%%%%%%%%%%%%%%%%%%%%%%%%%%%%%%%%%%%%%%%%%%%%%%%%%%%%
%%%
\begin{table}[htbp]
\begin{tabular}{|c|c|c|c|}
\hline 113,128: 66.9\% & 87,89: 25.9 \%& 53,55: 25.6 \% & 27,27: 23.6 \%\\
\hline 29,29:   22.9\% & 39,40: 21.1 \%& 27,28: 20.9 \% & 67,69: 16.6 \%\\
\hline 115,128: 16.6 \% & 36,41: 15.0 \%& 36,40: 13.0\%  & 26,28: 12.3\%\\
\hline  67,70: 10.9 \%& 55,56: 10.3 \% & &\\
\hline
\end{tabular}
\centering\label{T:ompA} \caption{The target sites $R[i,j]$ of
\emph{ompA} (interacts with \emph{MicA}) whose probabilities larger
than $10^{-1}$.}
\end{table}
%%%
%%%%%%%%%%%%%%%%%%%%%%%%%%%%%%%%%%%%%%%%%%%%%%%%%%%%%%%%%%%%%%%%%%%%%
%%%
%%%
%%%%%%%%%%%%%%%%%%%%%%%%%%%%%%%%%%%%%%%%%%%%%%%%%%%%%%%%%%%%%%%%%%%%%%%%%
%%%
\begin{table}[htbp]
\begin{tabular}{|c|c|c|c|}
\hline 52,60: 83.0\% & 15,17: 54.6 \%& 38,47: 24.7 \% & 15,16: 19.0 \%\\
\hline 72,75: 17.4\% & 77,78: 16.8 \%& 45,47: 14.2 \% & 71,74: 13.7 \%\\
\hline 73,75: 12.3 \% & 77,81: 11.1 \%& 14,17: 11.0\%  & \\
\hline
\end{tabular}
\centering\label{T:sodB} \caption{The target sites $R[i,j]$ of
\emph{sodB} (interacts with \emph{RyhB}) whose probabilities larger
than $10^{-1}$.}
\end{table}
%%%
%%%%%%%%%%%%%%%%%%%%%%%%%%%%%%%%%%%%%%%%%%%%%%%%%%%%%%%%%%%%%%%%%%%%%
%%%
%%%
%%%%%%%%%%%%%%%%%%%%%%%%%%%%%%%%%%%%%%%%%%%%%%%%%%%%%%%%%%%%%%%%%%%%%%%%%
%%%
\begin{table}[htbp]
\begin{tabular}{|c|c|c|c|}
\hline 87,93: 63.8\% & 39,48: 50.7 \%& 62,64: 43.6 \% & 70,72: 39.6 \%\\
\hline 30,30: 28.4\% & 70,73: 27.0 \%& 39,45: 17.0 \% & 87,92: 13.5 \%\\
\hline 40,45: 11.9 \% & 63,64: 11.4 \%&   & \\
\hline
\end{tabular}
\centering\label{T:fhlA} \caption{The target sites $R[i,j]$ of
\emph{fhlA} (interacts with \emph{OxyS}) whose probabilities larger
than $10^{-1}$.}
\end{table}
%%%
%%%%%%%%%%%%%%%%%%%%%%%%%%%%%%%%%%%%%%%%%%%%%%%%%%%%%%%%%%%%%%%%%%%%%
%%%
\subsection{Statistically generating interaction structure}
%%%
%%%%%%%%%%%%%%%%%%%%%%%%%%%%%%%%%%%%%%%%%%%%%%%%%%%%%%%%%%%%%%%%%%%%%
%%%
In this section, we generalize the idea of \cite{Ding:03} in order
to draw a representative sample from the Boltzmann equilibrium
distribution of RNA interaction structures. The section is divided
into two parts. At first we illustrate the correspondence between
the decomposition grammar, i.e.~the recursions for partition
functions and the sampling probabilities for mutually exclusive
cases and secondly, we describe the sampling algorithm.

The calculation of the sampling probabilities is based on the
recurrences of the partition functions since for mutually exclusive
and exhaustive cases, the key observation is that sampling
probability for a case is equivalent to the contribution to
partition function by the case divided by the partition function.
For instance, again we consider the decomposition of $J^{RT, {\sf
MK}}_{i,j;h,\ell}$. Set $\mathbb{P}^{0}_{i_1,j_1}$,
$\mathbb{P}^{1}_{i_1,j_1}$, $\mathbb{P}^{2}_{i_1,j_1}$,
$\mathbb{P}^{3}_{i_1,j_1}$ and $\mathbb{P}^{4}_{i_1,j_1}$ be the
sampling probabilities for all five cases showed in Figure
\ref{F:parexp3} ({\bf I}) anticlockwise respectively, then we have:
\begin{eqnarray*}
\mathbb{P}^{0}_{i_1,j_1}&=&Q^{\sf Hy,KK}_{i,i_1;h,h_1}Q^{RT, {\sf
KKA}}_{i_1+1,j;h_1+1,\ell}/Q^{RT, {\sf MK}}_{i,j;h,\ell},\\
\mathbb{P}^{1}_{i_1,j_1}&=&Q^{\sf Hy,KK}_{i,i_1;h,h_1}Q^{RT, {\sf
KF}}_{i_1+1,j;h_1+1,\ell}/Q^{RT, {\sf MK}}_{i,j;h,\ell},\\
\mathbb{P}^{2}_{i_1,j_1}&=&Q^{\sf Hy,KK}_{i,i_1;h,h_1}Q^{RT, {\sf
FF}}_{i_1+1,j;h_1+1,\ell}/Q^{RT, {\sf MK}}_{i,j;h,\ell},\\
\mathbb{P}^{3}_{i_1,j_1}&=&Q^{\sf Hy,KK}_{i,i_1;h,h_1}Q^{RT, {\sf
FK}}_{i_1+1,j;h_1+1,\ell}/Q^{RT, {\sf MK}}_{i,j;h,\ell},\\
\mathbb{P}^{4}_{i_1,j_1}&=&Q^{{\sf Hy KK}}_{i,j;h,\ell}/Q^{RT,{\sf
MK}}_{i,j;h,\ell}.
\end{eqnarray*}
Since the probabilities of all mutually exclusive and exhaustive
cases sum up to 1, we have
$\sum_{i_1,j_1}\mathbb{P}^{0}_{i_1,j_1}+\mathbb{P}^{1}_{i_1,j_1}+
\mathbb{P}^{2}_{i_1,j_1}+\mathbb{P}^{3}_{i_1,j_1}+
\mathbb{P}^{4}_{i_1,j_1}=1$, which coincides with
eqn.~(\ref{E:RTMK}).

Next we give a description of the sampling algorithm, as a
generalization of \cite{Ding:03}, we still take two stacks $A$ and
$B$. Stack $A$ stores sub-joint structures and their types $\xi$ in
the form of $\{(i,j;h,\ell; \xi)\}$, such as $(i,j;h,\ell; RT{\sf
MK})$ represents a sub-joint structure $J^{RT,{\sf
MK}}_{i,j;h,\ell}$. Stack $B$ collects interior/exterior arcs and
unpaired bases that will define a sampled interaction structure once
the sampling process finishes. At the beginning, $(1,N;1,M,
arbitrary)$ is the only element in stack $A$. A sampled interaction
structure is drawn recursively as follows: at first, start with
$(1,N;1,M; arbitrary)$, sample a pair of separated secondary
structures or a rts $(i,N;j,M; RT{\sf EE})$ according to their
sampling probabilities. In the former case, $(1,N;sec)$ and
$(1,M;sec)$ are stored in stack $A$. Otherwise, $(1,i-1;sec)$,
$(1,j-1;sec)$ and $(i,N;j,M; RT{\sf EE})$ are stored in stack $A$.
Secondly, given a new element in stack $A$, denoted by
$\{(i,j;h,\ell; \xi)\}$, we draw a particular case from all the
mutual exclusive and exhaustive cases according to the sampling
probabilities and store the corresponding sub-joint structures into
stack $A$, and all the interior arc, exterior arc or unpaired bases
sampled in the process will be stored in stack $B$. I.e.~after the
completion of sampling for a ``bigger" joint structure from stack
$A$ and storage of ``smaller" sub-joint structures derived in the
former process in stack $A$, also the storage of the sampled arcs
and unpaired bases of stack $B$, the element in the bottom of stack
$A$ is chosen to do the subsequent sampling. The who process
terminates when stack $A$ is empty and at the same time, a sampled
interaction structure formed in stack $B$.

%%%
%%%%%%%%%%%%%%%%%%%%%%%%%%%%%%%%%%%%%%%%%%%%%%%%%%%%%%%%%%%%%%%%%%%%%
%%%
\end{methods}

%%%
%%%%%%%%%%%%%%%%%%%%%%%%%%%%%%%%%%%%%%%%%%%%%%%%%%%%%%%%%%%%%%%%%%%%%
%%%
\section{Results and conclusions}
%%%
%%%%%%%%%%%%%%%%%%%%%%%%%%%%%%%%%%%%%%%%%%%%%%%%%%%%%%%%%%%%%%%%%%%%%
%%%
The complete set of recursions comprises for ts
$Q^{\bigtriangleup,\bigtriangledown,\square}_{i,j;r,s}$, 15
4D-arrays respectively, for right-tight structures
$Q^{RT}_{i,j;r,s}$, 20 4D-arrays, for dts $Q^{DT}_{i,j;r,s}$ and 20
4D-arrays. In addition, we need the usual matrices for the secondary
structures $R$ and $S$, and the above mentioned matrices for kissing
loops. The full set of recursions is compiled in the SM, Section~3.

%%%
%%%%%%%%%%%%%%%%%%%%%%%%%hybrid-matrix%%%%%%%%%%%%%%%%%%%%%%%%%%%%%%%%%%%%%%%%%%
%%%
\begin{figure}[ht]
\begin{center}
\includegraphics[width=\columnwidth]{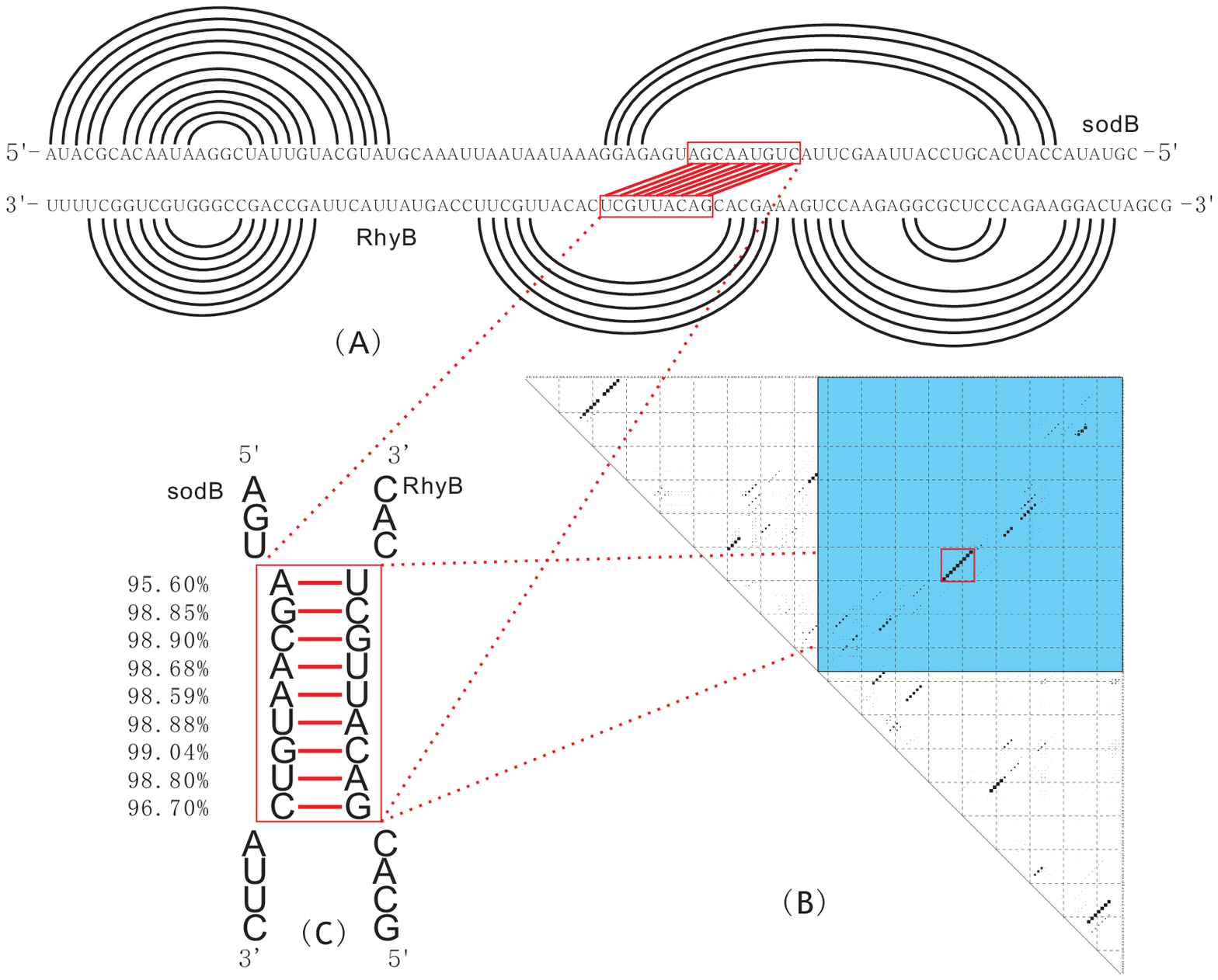}
\end{center}
\caption{\small (A) The natural  structure of
\emph{sodB}-\emph{RhyB} \citep{Geissmann}. (B) The
hybrid-probability matrix  generated via {\tt rip}2.0. This matrix
represents all potential contact regions of the \emph{sodB}
structure as squares, whose area is proportional to their respective
probability.  (C) ``Zoom'' into the most likely interaction region
as predicted by {\tt rip}2.0. All base pairs of the hybrid are
labeled by their probabilities.} \label{F:disscussion1}
\end{figure}
%%%
%%%%%%%%%%%%%%%%%%%%%%%%%%%%%%%%%%%%%%%%%%%%%%%%%%%%%%%%%%%%%%%%%%%
%%%

%%%
%%%%%%%%%%%%%%%%%%%%%%%%%hybrid-matrix%%%%%%%%%%%%%%%%%%%%%%%%%%%%%%%%%%%%%%%%%%
%%%
\begin{figure}[ht]
\begin{center}
\includegraphics[width=\columnwidth]{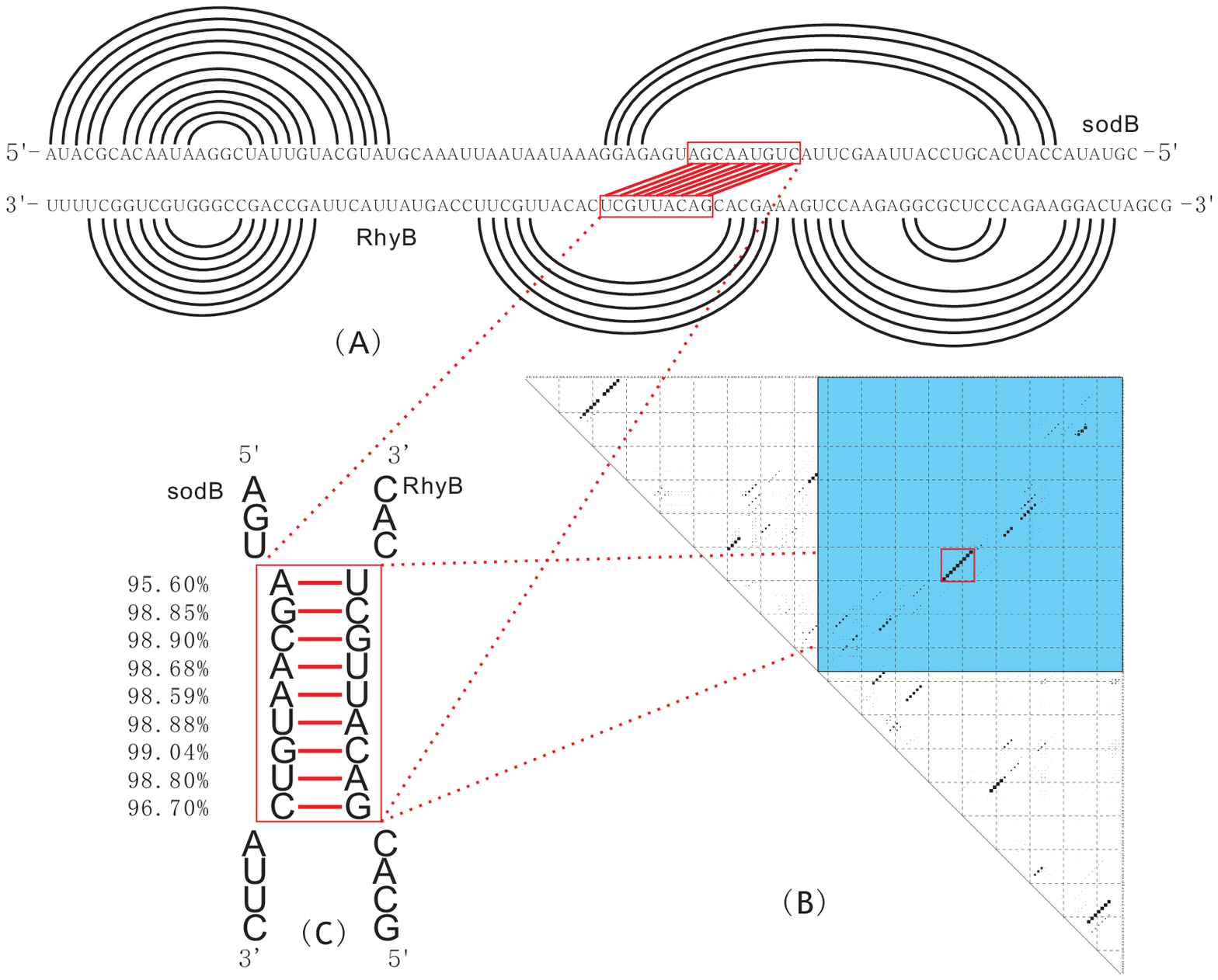}
\end{center}
\caption{\small (A) The natural  structure of
\emph{sodB}-\emph{RhyB} \citep{Geissmann}. (B) The
hybrid-probability matrix  generated via {\tt rip}2.0. This matrix
represents all potential contact regions of the \emph{sodB}
structure as squares, whose area is proportional to their respective
probability.  (C) ``Zoom'' into the most likely interaction region
as predicted by {\tt rip}2.0. All base pairs of the hybrid are
labeled by their probabilities.} \label{F:disscussion2}
\end{figure}
%%%
%%%%%%%%%%%%%%%%%%%%%%%%%%%%%%%%%%%%%%%%%%%%%%%%%%%%%%%%%%%%%%%%%%%
%%%

\begin{methods}
%%%
%%%%%%%%%%%%%%%%%%%%%%%%%%%%%%%%%%%%%%%%%%%%%%%%%%%%%%%%%%%%%%%%%%%%%
%%%
\section*{Acknowledgements}
%%%
%%%%%%%%%%%%%%%%%%%%%%%%%%%%%%%%%%%%%%%%%%%%%%%%%%%%%%%%%%%%%%%%%%%%%
%%%
We thank Bill Chen and Sven Findei{\ss} for comments on the
manuscript. This work was supported by the 973 Project of the
Ministry of Science and Technology, the PCSIRT Project of the
Ministry of Education, and the National Science Foundation of China
to CMR and his lab, grant No.\ STA 850/7-1 of the Deutsche
Forschungsgemeinschaft under the auspices of SPP-1258 ``Small
Regulatory RNAs in Prokaryotes'', as well as the European Community
FP-6 project SYNLET (Contract Number 043312) to PFS and his lab.
\end{methods}

\bibliographystyle{bioinformatics}
\bibliography{riphybrid}

\end{document}